\documentclass[12pt,reqno]{amsart}

\usepackage{amsrefs}
\usepackage{amssymb}
\usepackage{mathrsfs}
\usepackage{fancyhdr}
\usepackage{bbm}
\usepackage{comment}
\usepackage{xcolor}
\usepackage{hyperref}
\hypersetup{
  colorlinks=true,
  linkcolor=blue,
  citecolor=blue
}

\newtheorem{theorem}{Theorem}
\newtheorem*{theorem*}{Theorem}
\newtheorem{lemma}[theorem]{Lemma}
\newtheorem{proposition}[theorem]{Proposition}
\newtheorem{claim}[theorem]{Claim}

\newtheorem{maintheorem}{Theorem}

\theoremstyle{definition}

\newtheorem*{definition*}{Definition}
\newtheorem*{lemma*}{Lemma}
\usepackage{graphicx}
\usepackage{pstricks, enumitem, pst-node, pst-text, pst-plot}
\usepackage[capitalize]{cleveref}

\numberwithin{equation}{section}
\numberwithin{theorem}{section}

\newcommand{\R}{\mathbb{R}}
\newcommand{\N}{\mathbb{N}}

\newcommand{\eps}{\varepsilon}
\renewcommand{\P}{\mathcal{P}}	% probability measures

\usepackage{xparse}
\DeclareDocumentCommand\Pr{ m g }{\ensuremath{
    {   \IfNoValueTF {#2}
      {\mathbb{P}\left[{#1}\right]}
      {\mathbb{P}\left[{#1}\middle\vert{#2}\right]}%
    }
}}
\DeclareDocumentCommand\E{ m g }{\ensuremath{
    {   \IfNoValueTF {#2}
      {\mathbb{E}\left[{#1}\right]}
      {\mathbb{E}\left[{#1}\middle\vert{#2}\right]}%
    }
}}

\def\dd{\mathrm{d}}
\def\ee{\mathrm{e}}

\usepackage{todonotes}

\begin{document}

\title[]{Monotone additive statistics on heavy-tailed convolution semigroups}

\author{Tobias Fritz}
\address{University of Innsbruck, Innsbruck, Tyrol, Austria}
\email{tobias.fritz@uibk.ac.at}

\author{Xiaosheng Mu}
\address{Princeton University, Princeton, New Jersey, USA}
\email{xmu@princeton.edu}

\author{Omer Tamuz}
\address{California Institute of Technology, Pasadena, California, USA}
\email{tamuz@caltech.edu}


\thanks{\emph{Acknowledgments.} We thank two anonymous referees for their helpful comments and suggestions. This work was supported in part by NSF grant DMS-1944153 (OT).}

\date{\today}

\begin{abstract}
	We study sub-semigroups of the semigroup of probability measures on $\R$ and monotone additive statistics on them, by which we mean maps to the reals that are monotone with respect to the stochastic order and additive under convolution. We show that scalar multiples of the expectation are the unique monotone additive statistics on the semigroup of measures with finite $p$-th moment, for any $1 \le p < \infty$. We also prove that the entire semigroup of probability measures admits no non-zero monotone additive statistic at all.
\end{abstract}

\maketitle

\section{Introduction}
The set $\P(\R)$ of probability measures on the reals---which we just denote $\P$ for brevity---carries a natural semigroup structure given by the operation of convolution. It also admits a natural partial order: the stochastic order, also known as stochastic dominance. In this paper we study {\em monotone additive statistics}: maps from $\P$, or from sub-semigroups thereof, to $\R$, which respect the stochastic order and which are additive with respect to convolutions.

For partially ordered commutative semigroups (and similarly for other partially ordered algebraic structures like commutative rings), the set of order-preserving homomorphisms to $\R$ is a basic dual object whose study often yields fruitful insights~\cites{fritz, goodearl2010partially, krivine1964anneaux}. 
For $\P$ and its sub-semigroups in particular, studying such homomorphisms $\varphi$ is also naturally motivated in terms of applied probability. Namely, we can think of such $\varphi$ as a ``summary statistic''---a single number $\varphi(\mu)$ that captures some important property of a distribution $\mu \in \P$.
%The semigroup $\P$---considered without the stochastic order---is of course already an important object on its own.
The problem of finding well-behaved summary statistics arises in statistics, economics, operations research and other fields. For example, in financial asset pricing, $\mu$ can describe the distribution of returns of an asset. Then what price $\varphi(\mu)$ should we assign to the asset? If the mass of $\mu$ is below that of $\nu$ in the sense of first-order stochastic dominance, then we certainly expect $\varphi(\mu) \le \varphi(\nu)$. While if an asset is a portfolio consisting of two other assets $\mu$ and $\nu$ assumed independent, then its return distribution is described by the convolution $\mu \ast \nu$, and---under some assumptions---we would expect the prices to additionally satisfy $\varphi(\mu \ast \nu) = \varphi(\mu) + \varphi(\nu)$. 
We call $\varphi$ satisfying these properties \emph{monotone additive statistics}.
A similar approach is taken in the study of risk measures (see, e.g.,~\cite{goovaerts2004comonotonic}).

Indeed the homomorphisms $\varphi \colon \P \to \R$, and similarly for various sub-semigroups of $\P$, have been studied in the literature (see, e.g.,~\cites{ruzsa1988algebraic, mattner2004cumulants}).
However, despite the importance of monotonicity as described above, there seems to be less work on \emph{monotone} homomorphisms.
Through the present work we aim to fill part of this gap.

For $0 < p < \infty$, denote by $\P^p \subset \P$ the sub-semigroup of measures with finite $p$-th moment. If $p \ge 1$, then a monotone homomorphism $\P^p \to \R$ is given by the expectation $\E{\mu} = \int x\,\dd\mu(x)$. Some natural questions are: are there other monotone additive statistics on $\P^p$?
Are there \emph{any} monotone additive statistics on $\P^p$ for $p < 1$ or on all of $\P$ at all?
We answer these and similar questions in this paper: the expectation is (up to scalar multiples) the unique monotone additive statistic on $\P^p$ for $1 \le p < \infty$, and there are no non-zero monotone additive statistics on any $\P^p$ with $p < 1$ or on $\P$ at all.
Perhaps surprisingly, we show that there are semigroups strictly between $\P^1$ and $\P$ which \emph{do} admit non-trivial monotone additive statistics given by measuring the heaviness of tails. We end the paper with a number of further examples and open questions.

\subsection{Related literature}

Ruzsa and Sz{\'e}kely studied the semigroup $\P$ in their book ``Algebraic Probability Theory''  \cite{ruzsa1988algebraic}. They devote a chapter to additive statistics, where one of the questions they tackle is whether the expectation can be extended to a homomorphism on all of  $\P$. They show in their Theorem 2.4 that there are such extensions, and moreover there are extensions that assign 0 to all symmetric distributions on $\R$. However, there are no extensions that are non-negative for distributions supported on $\R_+$, and in particular there are no extensions that are monotone in our sense.

Mattner~\cite{mattner2004cumulants} studies the sub-semigroup $\bigcap_{p \in [1,\infty)} \P^p$ of measures that have all moments. He endows it with the topology of convergence in total variation, and of pointwise convergence of all moments. He shows that linear combinations of cumulants are the only continuous additive statistics on this semigroup to $\R$.

The same semigroup is studied in \cite{PST}, now as a partially ordered semigroup with respect to the stochastic order. It was shown there that for every $\mu,\nu \in \bigcap_p \P^p$ such that $\E{\mu} > \E{\nu}$, there is an $\eta \in \bigcap_{p \in [1,\infty)} \P^p$ such that $\mu \ast \eta \geq \nu \ast \eta$. This result immediately implies that if $\varphi$ is a monotone additive statistic and $\E{\mu} > \E{\nu}$, then $\varphi(\mu) \geq \varphi(\nu)$.
It follows that the scalar multiples of the expectation are the unique monotone additive statistics on $\bigcap_{p \in [1,\infty)} \P^p$.

The monotone additive statistics of the sub-semigroup of compactly supported measures is studied in \cite{mu2020monotone}. Let us denote it by $\P_c$. Then there are many monotone additive statistics on $\P_c$. For $t \in \R \setminus \{0\}$, Denote the normalized cumulant-generating function by 
$$
  \hat K_\mu(t) =  \frac{1}{t}\log\int\ee^{t x}\,\dd\mu(x),
$$  
and let $\hat K_\mu(0)$, $\hat K_\mu(-\infty)$ and $\hat K_\mu(\infty)$ respectively denote the expectation, the essential minimum, and the essential maximum of $\mu$. Then $t\mapsto \hat K_\mu(t)$ is continuous on $\R \cup \{-\infty,\infty\}$.

Any $t \in \R \cup\{-\infty,\infty\}$ determines
a monotone additive statistic on $\P_c$, given by $\mu \mapsto \hat K_\mu(t)$. It is shown in \cite{mu2020monotone}*{Theorem 1} that there are essentially no more, in the sense that the closed convex cone generated by these additive statistics coincides with the set of monotone additive statistics. Equivalently, for each monotone additive statistic $\varphi$ there is a finite Borel measure $\sigma$ on $\R \cup \{-\infty,\infty\}$ such that
\begin{align}
    \label{eq:avg-K}
  \varphi(\mu) = \int \hat K_\mu(t)\,\dd\sigma(t)
\end{align}
for all $\mu \in \P_c$.

For the semigroup of measures $\mu$ for which $\hat K_\mu(t)$ is finite for all $t \in \R$ (i.e., measures with a moment generating function), it is shown in \cite{mu2020monotone}*{Theorem 2} that every monotone homomorphism $\varphi$ with $\varphi(\delta_1)=1$ is likewise of the form \eqref{eq:avg-K}, with the additional restriction that $\sigma$ is compactly supported on $\R$.

\section{Main definitions and results}

We consider $\P$, the set of Borel probability measures on $\R$, as a partially ordered set with respect to the \emph{stochastic order}, also known as \emph{stochastic dominance}. Under this partial order, we have $\mu \le \nu$ if and only if their cumulative distributions functions are ordered pointwise,
\[
	\mu(  (-\infty,x] ) \ge \nu( (-\infty,x] ) \qquad \forall x \in \R.
\]
Equivalently, $\mu \leq \nu$ if there exists a standard probability space with random variables $X$ and $Y$ having distributions $\mu$ and $\nu$, and such that $X \leq Y$ almost surely. Equivalently again, $\mu \leq \nu$ if $\int f \,\dd\mu \leq \int f\,\dd\nu$ for all bounded non-decreasing $f \colon \R \to \R$. Intuitively, $\mu \leq \nu$ if one can arrive at $\mu$ by starting with $\nu$ and shifting mass to the left. If a measurable map $\pi \colon \R \to \R$ satisfies $\pi(x) \leq x$ for all $x$, then for any $\nu$ it holds that the push-forward $\pi_*\nu$ is dominated by $\nu$. 
      
Considering $\P$ as a semigroup under convolution and $\R$ as a semigroup under addition, we say that a map $\varphi \colon \P \to \R$ is an \emph{additive statistic} if it satisfies the homomorphism property $\varphi(\mu \ast \nu) = \varphi(\mu) + \mu(\nu)$. We say that it is {\em monotone} if $\mu \leq \nu$ implies $\varphi(\mu) \leq \varphi(\nu)$. For $x \in \R$ denote by $\delta_x$ the point mass at $x$. We then say that $\varphi \colon \P \to \R$ is {\em translation invariant} if $\varphi(\mu * \delta_x) = \varphi(\mu)$ for all $\mu \in \P$ and $x \in \R$. 

For $0 < p < \infty$, we denote by $\P^p \subset \P$ the sub-semigroup of measures $\mu$ that have finite $p$-th moment,
\[
	\int |x|^p\,\dd\mu(x)<\infty.
\]
For $p \ge 1$, Minkowski's inequality shows that this is indeed a sub-semigroup. For $p < 1$, this follows from the fact that $|x+y|^p \leq |x|^p + |y|^p$.

We also consider $\P_{\mathrm{Cram}}$, the semigroup of measures $\mu$ whose moment-generating function $M_\mu(t) = \int\ee^{tx}\,\dd\mu(x)$ is finite for $t$ in some neighborhood of zero which may depend on $\mu$ (Cram\'er's condition). This is a semigroup because of the multiplicativity of the moment-generating function under convolution. Note that it is a smaller semigroup than any $\P^p$, but larger than $\P_c$.

Our first main result shows in particular that the expectation is the unique monotone additive statistic on these semigroups, where for $\P^p$ we need to assume $p \ge 1$ in order for the expectation to exist.

\begin{maintheorem}
	\label{thm:L1}
	On the following sub-semigroups of $\P$, the monotone additive statistics $\varphi$ are precisely the maps of the form $\varphi(\mu) = c \E{\mu}$ for some $c \geq 0$:
	\begin{enumerate}[label=(\roman*)]
		\item\label{item:Lpgt1} $\P^p$ for $1 \le p < \infty$.
		\item\label{item:cramer} $\P_{\mathrm{Cram}}$. 
	\end{enumerate}
\end{maintheorem}

Our second main result is a non-existence result for monotone additive statistics on two larger types of semigroups.

\begin{maintheorem}
	\label{thm:PR}
	On the following sub-semigroups of $\P$, the only monotone additive statistic $\varphi$ is $\varphi = 0$:
	\begin{enumerate}[label=(\roman*)]
		\item\label{item:PR} $\P$ itself.
		\item\label{item:Lplt1} $\P^p$ for $0 < p < 1$.
	\end{enumerate}
\end{maintheorem}

We also show that there are intermediate semigroups, situated between $\P^1$ and $\P$, which \emph{do} admit monotone additive statistics (necessarily different from the expectation).
To this end, let $a = (a_n)_{n \in \N}$ be a sequence of positive reals.
Consider the set $\P^a \subseteq \P$ of measures $\mu \in \P$ for which the limit
\begin{align*}
	\psi_a(\mu) := \lim_{n \to \infty} a_n \cdot \mu((n,\infty))
\end{align*}
exists and is finite. A good example is $a_n=n$, in which case $\psi_a(\mu) = c$ if and only if the right tail of $\mu$ decays like $c/n$.
We think of this $\psi_a$ as capturing the ``thickness of the tail'', as compared to $1/n$, and similarly for other choices of $a$.

Our main interest will be in sequences $(a_n)$ that are non-decreasing and grow sufficiently mildly. As the next theorem shows, under these conditions, $\psi_a$ is a monotone additive statistic. If the sequence grows slowly enough, then $\P^a$ includes all distributions with finite expectation.
\begin{maintheorem}
\label{thm:Lpsi}
    Suppose that $(a_n)_n$ is non-decreasing and
    \begin{enumerate}
        \item $\lim_n \frac{a_n}{a_{n+1}} = 1$, and 
	\item $\sup_n \frac{a_{2n}}{a_n} < \infty$. 
    \end{enumerate}
    Then $\P^a$ is a semigroup and $\psi_a : \P^a \to \R$ is a monotone additive statistic.
    If additionally $\sup_n \frac{a_n}{n} <\infty$, then $\P^1 \subsetneq \P^a$.
\end{maintheorem}

We also investigate a few other examples of intermediate semigroups, with diverse behavior of the sets of monotone additive statistics. For example, we show that there exists a semigroup in which the set of monotone additive statistics contains an infinite-dimensional vector space corresponding to the solutions of the Cauchy functional equation.

\section{Proof of Theorem~\ref{thm:L1}: Uniqueness of the expectation}

Suppose $\mu$ and $\nu$ are not comparable in the stochastic order. Is it possible that there exists a measure $\eta$ such that $\mu \ast \eta \geq \nu \ast \eta$? This question was first asked, independently, by Tarsney~\cite{tarsney2018exceeding} and Pomatto et al.~\cite{PST}. The latter show that this happens whenever $\E{\mu} > \E{\nu}$, using a lemma due to Ruzsa and Sz{\'e}kely~\cite{ruzsa1988algebraic}.

In general, it is impossible to control $\eta$ without controlling $\mu$ and $\nu$. In particular, even when $\mu,\nu \in \P^1$, there may not exist such an $\eta$ in $\P^1$. The next lemma, which is a key ingredient of our proof of Theorem~\ref{thm:L1}, shows that when $\mu,\nu \in \P_{\mathrm{Cram}}$, then one can take $\eta$ equal to a Laplace distribution, and in particular to also be in $\P_{\mathrm{Cram}}$. This was shown for compactly supported $\mu,\nu$ in \cite{mu2020background}*{Theorem 1} and \cite{PST}*{Theorem 3}.

For $r > 0$, let the \emph{Laplace measure with parameter $r$} be the probability measure on the real line having density function $h(x)=\frac{r}{2}\ee^{-r|x|}$.
\begin{lemma}
\label{lem:catalytic}
Suppose that $\mu,\nu \in \P_{\mathrm{Cram}}$ satisfy $\E{\mu} > \E{\nu}$. Let $\eta_r$ be the Laplace measure with parameter $r$. Then
\begin{align*}
    \mu \ast \eta_r \geq \nu \ast \eta_r
\end{align*}
for every small enough $r > 0$.
\end{lemma}

\begin{proof}
  Denote the c.d.f.s of $\mu$ and $\nu$ by  $F(x) = \mu((-\infty,x])$ and $G(x) = \nu((-\infty,x])$, respectively.
  
  Since $\mu,\nu \in \P_{\mathrm{Cram}}$, there is some $s > 0$ such that the moment-generating functions $M_\mu$ and $M_\nu$ are finite on $[-s,s]$. It follows from integration by parts that for $r \in [-s,s]$, we have
  \begin{align*}
    M_\mu(r) - M_\nu(r) = r\int_{-\infty}^\infty \ee^{r x}[G(x)- F(x)]\,\dd x.
  \end{align*}
  Hence
  \begin{equation}
    \label{mfg_dominating}
    \int_{-\infty}^\infty \ee^{s |x|} [G(x) - F(x)]\,\dd x
  \end{equation}
  is well defined and finite.

  For $r \in [0,s]$ and $y \in \R$, let $e_{y,r} \colon \R \to \R$ be given by
  \begin{align*}
    e_{y,r}(x) = \ee^{r |y|}\ee^{-r|y-x|}.
  \end{align*}
  We make two observations:
  \begin{enumerate}
  \item $0 \leq e_{y,r}(x) \leq \ee^{s |x|}$ for any $r \in [0,s]$ and $y \in \R$.
  \item Suppose that $(r_n)_n$ is a sequence in $[0,s]$ that converges to $0$, and $(y_n)_n$ is any sequence. Then the sequence of functions $(e_{y_n,r_n})_n$ converges pointwise to $1$.
  \end{enumerate}
  Given these observations, it follows from dominated convergence\footnote{In more detail, the signs of the integrands in the sequence are the same for every $x$, and given by the sign of $G(x) - F(x)$. One can therefore apply dominated convergence to the positive and negative parts separately, with the corresponding parts of~\eqref{mfg_dominating} as the dominating integral.} that for any sequences $r_n \to 0$ and $y_n$,
  \begin{align}
    \label{eq:lim-r-positive}
    \lefteqn{\lim_n \int_{-\infty}^\infty e_{y_n,r_n}(x)[G(x)-F(x)]\,\dd x} \nonumber\\
    &=
    \int_{-\infty}^\infty [G(x)-F(x)]\,\dd x = \E{\mu}-\E{\nu} > 0.
  \end{align}

  We claim that for all $r$ small enough it holds that
  \begin{align*}
    \Delta_r(y) := \int_{-\infty}^\infty \ee^{-r|y-x|}[G(x)-F(x)]\,\dd x
  \end{align*}
  is positive for all $y$. If not, then for any sequence $r_n$ tending to zero we can choose a sequence $y_n$ so that $\Delta_{r_n}(y_n) \leq 0$ for all $n$. It follows that $\ee^{r_n |y_n|}\Delta_{r_n}(y_n) \leq 0$, which contradicts \eqref{eq:lim-r-positive}.

  Fix some such $r$ small enough, and let $h(x) = \frac{r}{2}\ee^{-r |x|}$ be the p.d.f.\ of the Laplace measure $\eta$ with parameter $r$. The c.d.f.\ of $\mu \ast \eta$ is $h \ast F$, and the c.d.f.\ of $\nu \ast \eta$ is $h \ast G$. Furthermore
  \begin{align*}
    h \ast G - h \ast F = \frac{r}{2}\Delta_r > 0,
  \end{align*}
  and so $\mu \ast \eta \geq \nu \ast \eta$.
\end{proof}

\begin{lemma}
\label{lem:point-masses}
Let $\mathcal{S} \subset \P$ be a sub-semigroup that contains the point masses $\{\delta_c \,:\, c \in \R\}$, and let $\varphi \colon \mathcal{S} \to \R$ be a monotone additive statistic. Then $\varphi(\delta_c) = c \varphi(\delta_1)$.
\end{lemma}
\begin{proof}
By additivity, $\varphi(\delta_c) = c \varphi(1)$ for rational $c$. By monotonicity it follows that this holds for all $c$.
\end{proof}

\begin{lemma}
\label{lem:mixtures}
Let $\varphi \colon \P^p \to \R$ be a monotone additive statistic, and let $\mu,\nu \in \P_{c}$ have the same expectation. Then for every $\zeta \in \P^p$ and $\alpha \in [0,1]$, it holds that 
\[
\varphi(\alpha \mu + (1-\alpha) \zeta) = \varphi(\alpha \nu + (1-\alpha) \zeta). 
\]
In particular, $\varphi(\mu)=\varphi(\nu)$, and so $\varphi(\mu) = \lambda \E{\mu}$ for some $\lambda \geq 0$.
\end{lemma}
\begin{proof}
It suffices to prove the inequality $\ge$.
Fix some $c >0$. Then  $\E{\mu \ast \delta_c} > \E{\nu}$, and so by Lemma~\ref{lem:catalytic} there exists a Laplace measure $\eta \in \P^p$ such that $(\mu \ast \delta_c) \ast \eta \geq \nu \ast \eta$. Clearly we also have $(\zeta \ast \delta_c) \ast \eta \geq \zeta \ast \eta$. So 
\[
\alpha (\mu \ast \delta_c \ast \eta) + (1-\alpha)(\zeta \ast \delta_c \ast \eta) \geq \alpha (\nu \ast \eta) + (1-\alpha)(\zeta \ast \eta).
\]
Hence
\[
\big[\alpha \mu + (1-\alpha)\zeta\big]\ast \delta_c \ast \eta  \geq \big[\alpha \nu  + (1-\alpha)\zeta\big] \ast \eta.
\]
It now follows by the monotonicity and additivity  of $\varphi$ that 
\[
\varphi\big(\alpha \mu + (1-\alpha)\zeta\big)+ \varphi(\delta_c)  \geq \varphi\big(\alpha \nu  + (1-\alpha)\zeta\big).
\]
By Lemma~\ref{lem:point-masses} we have $\varphi(\delta_c) = c \varphi(\delta_1)$, and so letting $c \to 0$ yields the desired inequality.
\end{proof}

By Lemma~\ref{lem:catalytic}, part~\ref{item:cramer} of Theorem~\ref{thm:L1} is a consequence of the following observation.

\begin{proposition}
	\label{prop:catalytic}
	Let $\mathcal{S} \subseteq \P^1$ be any sub-semigroup such that:
	\begin{enumerate}
		\item $\mathcal{S}$ contains the point masses $\{\delta_c \,:\, c \in \R\}$.
		\item If $\mu,\nu \in \mathcal{S}$ satisfy $\E{\mu} > \E{\nu}$, then there is $\eta \in \mathcal{S}$ with $\mu \ast \eta \ge \nu \ast \eta$.
	\end{enumerate}
	Then the monotone additive statistics on $\mathcal{S}$ are precisely the scalar multiples of $\mathbb{E}$.
\end{proposition}

\begin{proof}
	Given a monotone additive statistic $\varphi \colon \mathcal{S} \to \R$, we show that $\varphi(\mu) = \E{\mu} \, \varphi(\delta_1)$ for every $\mu \in \mathcal{S}$.

	Fix arbitary $\eps > 0$. Then by assumption we have $\eta_{\pm} \in \mathcal{S}$ such that
	\[
		\delta_{\E{\mu} - \eps} \ast \eta_- \le \mu \ast \eta_-, \qquad
		\mu \ast \eta_+ \le \delta_{\E{\mu} + \eps} \ast \eta_+,
	\]
	so that applying $\varphi$ results in
	\[
		\varphi(\delta_{\E{\mu} - \eps}) \le \varphi(\mu) \le \varphi(\delta_{\E{\mu} + \eps}),
	\]
	or equivalently by Lemma~\ref{lem:point-masses},
	\[
		(\E{\mu} - \eps) \varphi(\delta_1) \le \varphi(\mu) \le (\E{\mu} + \eps) \varphi(\delta_1).
	\]
	This proves the claim in the limit $\eps \to 0$.
\end{proof}

We now focus on the more involved proof of part~\ref{item:Lpgt1} of Theorem~\ref{thm:L1}. Let $\varphi \colon \P^p \to \R$ be a monotone additive statistic. Write $\lambda = \varphi(\delta_1)$. By monotonicity and since $\varphi(0)=0$, we have $\lambda \geq 0$. From the above Lemma~\ref{lem:mixtures}, we deduce that for any compactly supported $\mu$,
\[
\varphi(\mu) = \varphi(\delta_{\E{\mu}}) = \lambda \E{\mu}.
\]
We next show that the same holds for any $\mu \in \P^p$ that is \emph{bounded from below} (i.e., $\mu([-M,\infty))=1$ for some $M$) but not bounded from above. Suppose for the sake of contradiction that there exists such a $\mu$ with $\varphi(\mu) \neq \lambda \E{\mu}$. Shifting $\mu$ by a constant if necessary, we can assume that $\mu$ is supported on $\R_{+}$, and $\varphi(\mu) \neq \lambda \E{\mu}$. For each positive integer $n$, consider the function $f^n$ with values 
\begin{align}
\label{eq:min}
    f^n(x) := \min\{x,n\}.
\end{align}
Since $f^n(x) \leq x$ for all $x \in \R$, the measure $\mu_n = f^n_*(\mu)$ satisfies $\mu_n \leq \mu$.

It follows from the Monotone Convergence Theorem that $\E{\mu_n} \to \E{\mu}$. By monotonicity, we have $\varphi(\mu) \geq \varphi(\mu_n) = \lambda \E{\mu_n}$ for each $n$. Hence $\varphi(\mu) \geq \lambda \E{\mu}$. Since by assumption equality does not hold, we deduce $\varphi(\mu) > \lambda \E{\mu}$. We can then choose a large $k$ such that $\varphi(\mu) \geq \lambda \E{\mu} + \frac{1}{k}$. Let $\nu_n = \mu^{(n k)}$ equal the convolution of $\mu$ with itself $n k$ times. Then additivity of $\varphi$ and $\mathbb{E}$ implies
\[
\varphi(\nu_n) \geq \lambda \E{\nu_n} + n. 
\]

For each $n$, we now choose positive number $a_n$ that is sufficiently large to satisfy the following properties:
\begin{enumerate}
    \item $p_n := \nu_n([0, a_n]) \geq 1 - \frac{1}{n}$;
    \item $\int_{a_n}^\infty x^p \, \dd \nu_n(x) \leq 2^{-n}$.
\end{enumerate}
Denote by $\nu_n^1$ and $\nu_n^2$ the measure $\nu_n$, conditioned on $[0,a_n]$ and $(a_n,\infty)$, respectively:
\begin{align*}
	\nu_n^1(A) &= \frac{\nu_n(A \cap [0,a_n])}{p_n} = \nu_n(A | [0,a_n]),\\[3pt]
	\nu_n^2(A) &= \frac{\nu_n(A \cap (a_n,\infty))}{1-p_n} = \nu_n(A | (a_n,\infty)).\\
\end{align*}
Then 
$$
  \nu_n = p_n \nu_n^1 + (1-p_n)\nu_n^2.
$$
Denote $m_n = \E{\nu_n^1}$, and let
$$
  \zeta_n = p_n \delta_{m_n} + (1-p_n)\nu^2_n.
$$
That is, $\zeta_n$ is obtained from $\nu$ by contracting all the mass below $a_n$ to an atom at its (conditional) expectation $m_n = \E{\nu_n^1}$. It thus follows from Lemma~\ref{lem:mixtures} that 
\begin{align*}
    \varphi(\zeta_n)  =\varphi(\nu_n).
\end{align*}

Next, let
\begin{align*}
  \eta_n := \zeta_n \ast \delta_{-m_n} = p_n\delta_0 + (1-p_n)\nu_n^2 \ast \delta_{-m_n}
\end{align*}
be the measure $\zeta_n$ shifted to the left by the conditional expectation; note that it is supported on $\R_{+}$ and is unbounded from above. By additivity, 
\begin{align*}
\varphi(\eta_n) &= \varphi(\zeta_n) - \lambda \E{\nu_n^1} \\
&\geq \varphi(\nu_n) - \lambda \E{\nu_n} \\
&\geq n,
\end{align*}
so that $\varphi(\eta_n)$ is large. On the other hand, its $p$-th moment is small: since
\begin{align*}
    \eta_n = p_n \delta_0 + (1-p_n)\nu_n^2 \ast \delta_{-m_n},
\end{align*}
it follows from the second property of $a_n$ that
\begin{align*}
    \int_0^\infty x^p\,\dd\eta_n(x) = \int_{a_n}^\infty (x-m_n)^p\,\dd\nu_n(x) <  \int_{a_n}^\infty x^p\,\dd\nu_n(x) \leq 2^{-n}.
\end{align*}
Moreover, the first condition on $a_n$ ensures that $\eta_n$ has a mass point of size at least $1-\frac{1}{n}$ at zero. 

Now define $F_n$ to be the c.d.f.\ of $\eta_n$, and consider $F(x) = \inf_{n} F_n(x)$.
Then clearly $F(x) = 0$ for $x < 0$.
For given $b \geq 0$, fix $N$ large enough so that $F_1(b) \le 1 - \frac{1}{N}$.
Then the defining infimum of $F(x)$ for $x \in [0,b]$ is achieved at some $n \in \{1,\dots,N\}$, since for $n > N$ we have
\[
	F_n(x) \geq F_n(0) \geq 1 - \frac{1}{n} > F_1(b) \ge F_1(x). 
\]
Hence each $F_n$ being non-decreasing and right-continuous yields the same properties for $F$. Moreover, the fact that $F_n(0) \geq 1-\frac{1}{n}$ and each $F_n(x) \to 1$ as $x \to \infty$ implies that $F(x) \to 1$ as $x \to \infty$. Hence $F$ is the c.d.f.\ of some probability measure $\eta$ supported on $\R_{+}$. In addition, $\eta \in \P^p$ because
\begin{align*}
\int x^p \,\dd\eta(x) &= \int_{0}^{\infty} px^{p-1}(1-F(x)) \,\dd x \\
&\leq \int_{0}^{\infty} px^{p-1}\left(\sum_{n = 1}^{\infty} (1 - F_n(x))\right) \,\dd x \\
&= \sum_{n = 1}^{\infty} \int_{0}^{\infty} px^{p-1}(1 - F_n(x)) \,\dd x \\
&= \sum_{n = 1}^{\infty} \int_0^\infty x^p\,\dd\eta_n(x) \\
&\leq \sum_{n = 1}^{\infty} 2^{-n},
\end{align*}
where we have used the fact that
\begin{align*}
    1 - F(x) = \sup_{n\geq 1} \, (1 - F_n(x)) \leq \sum_{n\geq 1}(1 - F_n(x)),
\end{align*}
as well as the Monotone Convergence Theorem. 

Now observe that since $F(x) \leq F_n(x)$ for all $x$, we have $\eta \geq \eta_n$. Therefore $\varphi(\eta) \geq \varphi(\eta_n) \geq n$ for each $n$. This contradicts the assumption that $\varphi(\eta) \in \R$. 

Therefore, for any $\mu$ that is bounded from below, we must have $\varphi(\mu) = \lambda \E{\mu}$. A symmetric argument shows that the same is true for any $\mu$ that is bounded from above. Finally, for a general $\mu$ that may be unbounded on both sides, we have by monotonicity, with $f^n(x) = \min\{x,n\}$ from \eqref{eq:min},
\[
\varphi(\mu) \geq \varphi(f^n_*\mu) = \lambda \E{f^n_*\mu},
\]
so that $\varphi(\mu) \geq \lambda \E{\mu}$ by letting $n \to \infty$ and applying dominated convergence. Likewise, if we denote $g^n(x) = \max\{x,-n\}$, then 
\[
\varphi(\mu) \leq \varphi(g^n_*\mu) = \lambda \E{g^n_*\mu},
\]
so that $\varphi(\mu) \leq \lambda \E{\mu}$ also holds. This concludes the proof of Theorem~\ref{thm:L1}.

\bigskip

As a final note, in \cite{PST} it was shown that the only monotone additive statistics on $\bigcap_{p\geq 1} \P^p$ are again the scalar multiples of expectation. 
The argument in the above proof can be adapted to reproduce that result.
For this we just need to choose $a_n$ sufficiently large so that 
\[
\int_{a_n}^\infty x^p\,\dd\nu_n(x) \leq 2^{-n} ~~ \text{for all} ~~ p \in [1, n].
\]
Such $a_n$ exists because by H{\"o}lder's inequality, it is enough to require this bound in the cases $p = 1$ and $p = n$.
Then $\int_0^\infty x^p \,\dd\eta_n(x)\leq 2^{-n}$ for  $p \leq n$ and for any $p$,
\[
\int x^p \,\dd\eta(x) \leq \sum_{n=1}^{\infty} \int x^p \,\dd\eta_n(x) < \infty,
\]
so that $\eta \in \bigcap_{p \geq 1} \P^p$ as well. %This proof shows that the Rusza-Szekely lemma is not essential for the characterization of monotone additive statistics on this space. 
%\end{comment}

\section{Proof of Theorem~\ref{thm:PR}}

In this subsection we prove Theorem~\ref{thm:PR} on the non-existence of monotone additive statistics on $\P$ and $\P^p$ for $0 < p < 1$. We treat the case of $\P$ itself first, again starting with a lemma.

As already mentioned, we say that a statistic $\varphi \colon \P \to \R$ is \emph{translation invariant} if $\varphi(\mu \ast \delta_c) = \varphi(\mu)$ for all $\mu \in \P$ and $c \in \R$. Note that such $\varphi$ does not have to be additive.
\begin{lemma}
  \label{lem:bounded}
	Let $\varphi : \P \to \R$ be any translation invariant and monotone map. Then $\varphi$ is bounded.
\end{lemma}
\begin{proof}
  Let $g \colon \R \to \R$ be given by $g(x) = \max\{0,x\}$. Given $\mu \in \P$, denote by $g_*\mu$ the push-forward of $\mu$ under $g$. Then $g_*\mu$ is supported on $[0,\infty)$. It is immediate that $\mu \leq g_*\mu$.
  
  By the translation invariance and monotonicity assumptions on $\varphi$, for all $\mu \in \P$ and $c \in \R$ it holds that $\varphi(\mu \ast \delta_c) = \varphi(\mu)$  and $\varphi(\mu) \leq \varphi(g_*\mu)$.
  
  Assume now by contradiction that $\varphi$ is unbounded from above, so that there is a sequence $(\mu_n)_{n\in\N}$ with $\varphi(\mu_n) \geq n$. Define the sequence of measures $(\nu_n)_n$ as follows: for each $n$, choose $a_n$ large enough so that $\mu_n((-\infty,a_n]) \geq 1-1/n$, and let $\nu_n = g_*(\mu_n \ast \delta_{-a_n})$ be the translation of $\mu_n$ by $-a_n$, pushed forward by $g$. Note that (i) $\nu_n([0,\infty))=1$, (ii) $\nu_n(\{0\}) \geq 1-1/n$, and (iii) $\varphi(\nu_n) \geq \varphi(\mu_n) \ge n$.

  Denote by
  \begin{align*}
    F_n(x) = \nu_n((-\infty,x]) = \nu_n([0,x])
  \end{align*}
  the c.d.f.\ of $\nu_n$, and let $F(x) = \inf_n F_n(x)$. As in the proof of Theorem~\ref{thm:L1}, $F$ is non-decreasing, right-continuous, and satisfies $F(x) \to 1$ as $x \to \infty$. Thus $F$ is the cumulative distribution function of some $\nu \in \P$. Since $F(x) \leq F_n(x)$ for all $x$ and $n$, we have that $\nu \geq \nu_n$ for all $n$, and so $\varphi(\nu) \geq \varphi(\nu_n) \geq n$ for all $n$. We have thus reached a contradiction.
  
  An analogous argument with respect to going downwards in the stochastic order shows that $\varphi$ must also be bounded below.
\end{proof}

Lemma~\ref{lem:bounded}  is clearly not true without the assumption of translation invariance: for example, for any $p \in (0,1)$, taking the quantile
\[
	\mu \longmapsto \inf \{ x \in \R \mid \mu( (-\infty,x] ) \ge p \}
\]
defines an unbounded monotone map $\P \to \R$.

We now return to the proof of Theorem~\ref{thm:PR}\ref{item:PR}. Since $\R$ has no nonzero torsion elements,
$\varphi$ can only be bounded if it is identically zero. Thus the claim will follow from Lemma~\ref{lem:bounded} if we can show that every monotone additive statistic $\varphi \colon \P \to \R$ is translation invariant. To this end, we need to show that for every $x \in \R$ it holds that $\varphi(\delta_x)=0$. By Lemma \ref{lem:point-masses} it suffices to show $\varphi(\delta_1)=0$.

Suppose $\varphi(\delta_1) > 0$. As in the proof of Theorem~\ref{thm:L1}, let $f^n(x) = \min\{x,n\}$, so that $f^n_*\mu \leq \mu$. Let $\mu$ be any measure that has infinite expectation and is supported on $[0,\infty)$. Then $\mu_n = f^n_*\mu$ has an expectation (since it has compact support), and $\lim_n \E{\mu_n} = \infty$ by monotone convergence. By Theorem~\ref{thm:L1}, $\varphi(\mu_n)=\varphi(\delta_1)\E{\mu_n}$, and so $\lim_n\varphi(\mu_n) =\infty$. But $\mu \geq \mu_n$, and so  $\varphi(\mu) \geq \varphi(\mu_n)$. This gives a contradiction and concludes the proof of part~\ref{item:PR} of Theorem~\ref{thm:PR}.

\medskip

Part~\ref{item:Lplt1} claims that there is likewise no non-trivial monotone additive statistic on $\P^p$ for $0 < p < 1$, which we prove now.
Such $\P^p$ contains a measure $\mu$ supported on $\R_+$  with infinite expectation. Approximating $\mu$ from below by $f^n_*\mu$, we deduce by monotonicity and the case of $\P^1$ from Theorem~\ref{thm:L1} that $\varphi(\mu) \geq \varphi(\delta_1) \, \E{f^n_* \mu}$ for all $n$. Since $\E{f^n_*\mu} \to \infty$ as $n \to \infty$, we obtain $\varphi(\delta_1) = 0$, and therefore $\varphi$ is translation invariant as above.

The rest of the argument proceeds similarly to the proof of Theorem~\ref{thm:L1}: if $\varphi(\mu) > 0$ for some positively supported $\mu \in \P^p$, then we can find a sequence $\nu_n \in \P^p$ such that $\varphi(\nu_n) \geq n$.
By choosing $a_n$ sufficiently large, we can make the resulting measures $\eta_n$ satisfy 
$\eta_n(\{0\}) \geq 1 - \frac{1}{n}$ and $\int x^p\,\dd\eta_n(x) \leq 2^{-n}$. Moreover, $\varphi(\eta_n) = \varphi(\nu_n) \geq n$ by translation invariance. Hence the stochastically dominant $\eta$
satisfies $\varphi(\eta) \geq \varphi(\eta_n) \geq n$ for each $n$, leading to a contradiction. Therefore $\varphi(\mu) = 0$ for any positively supported $\mu$, and thus for any $\mu \in \P^p$ that is bounded from below. Symmetrically $\varphi(\mu) = 0$ for any $\mu$ bounded from above. The general case follows from the same approximation argument as at the end of the proof of Theorem~\ref{thm:L1}.

\section{Proof of Theorem~\ref{thm:Lpsi}}.

	We first prove that $\sup_n \frac{a_n}{n} <\infty$ implies the strict containment $\P^1 \subset \P^a$.
	So suppose $\mathbb{E}[\mu] < \infty$, and write $s_n = \mu((n,\infty))$. Then
    \begin{align*}
        \sum_{n=1}^\infty s_n= \sum_{n=1}^\infty \mu((n,\infty)) \leq  \int_0^\infty \mu((x,\infty))\,\dd x < \infty.
    \end{align*}
    Since $s_n$ is decreasing and summable, it passes the Olivier Test, which means $\lim_n n s_n=0$.
    By $\sup_n \frac{a_n}{n} < \infty$, it follows that 
    \begin{align*}
        \psi_a(\mu)=\lim_n a_n s_n = \lim_n \frac{a_n}{n} n s_n = 0.
    \end{align*}
    We have thus shown that $\P^1 \subseteq \P^a$.

    For strict containment we consider two cases. First, if $a_n$ is bounded from above, then $\psi_a=0$ on all of $\P$, and thus $\P^a = \P \supsetneq \P_1$. 
    If $a_n$ is unbounded, then let $\mu$ be any probability measure with cumulative distribution function $F$ satisfying $F(n) = 1-\min\{1,a_n^{-1}\}$ for all $n \in \N$, so that $\psi_a(\mu)=1$. Thus $\mu \in \P^a$, and it remains to be shown that $\mu$ has infinite first moment. Note that
    \begin{align*}
        \int |x| \,\dd\mu(x) \geq \sum_{n=1}^\infty (1-F(n)) = \sum_{n=1}^\infty \min\{1,a_n^{-1}\}.
    \end{align*}
    Since $\sup_n \frac{a_n}{n} < \infty$, there is $\eps>0$ such that $a_n \leq \eps^{-1} n$ for all $n \ge 1$. Thus
    \begin{align*}
	    \int |x| \,\dd\mu(x) \geq \sum_{n=1}^\infty \min\{1,\eps n^{-1}\} = \infty.
    \end{align*}
    
  Concerning the properties of $\psi_a$, its monotonicity is immediate whenever it exists and is finite.
  We thus prove the claim by showing that $\psi_a(\eta) = \psi_a(\mu) + \psi_a(\nu)$ whenever $\eta = \mu \ast \nu$. We do this by estimating $a_n \cdot \eta((n,\infty))$ from both sides.

For this proof it will be useful to use probabilistic notation. Let $X$ and $Y$ be two independent random variables, with distributions $\mu$ and $\nu$. Hence their sum $X+Y$ has distribution $\eta$.

Fix any $\epsilon > 0$, and let $m$ be sufficiently large such that $\Pr{\vert X \vert \leq m}$ and $\Pr{\vert Y\vert \leq m}$ are both larger than $1-\epsilon$. Then for every $n \ge 0$,
\begin{align*}
  \eta((n,\infty)) 
  &= \Pr{X+Y > n}\\
  &\geq \Pr{\vert X \vert \leq m, ~Y > n + m} + \Pr{\vert Y \vert \leq m, ~X > n + m}  \\
&= \Pr{\vert X \vert \leq m} \cdot \Pr{Y > n + m} + \Pr{\vert Y \vert \leq m} \cdot \Pr{X > n + m} \\
&\geq (1-\epsilon) \left(\Pr{Y > n + m} + \Pr{X > n + m}\right).
\end{align*}
Thus 
\begin{align*}
 \lefteqn{a_n \cdot \eta((n,\infty)) }\\
 &\geq \frac{(1-\epsilon)a_n}{a_{n+m}} \cdot \left(a_{n+m}\Pr{Y > n + m} + a_{n+m}\Pr{X > n + m} \right)   
\end{align*}
By assumption $\lim_n \frac{a_n}{a_{n+1}} = 1$, which implies that $\lim_n \frac{a_n}{a_{n+m}} = 1$. Hence, 
letting $n \to \infty$ yields
\[
	\liminf_{n \to \infty} a_n \cdot \eta((n,\infty)) \geq (1-\epsilon)(\psi_a(\mu)+\psi_a(\nu)).
\]
Since $\epsilon$ is arbitrary, we have $\liminf_{n \to \infty} a_n \cdot \eta((n,\infty)) \geq \psi_a(\mu)+\psi_a(\nu)$.

In the opposite direction, for $n > 2m$ we can write 
\begin{align*}
	\Pr{X + Y > n} ={} & \Pr{X \leq m, ~X + Y > n} \\
			& + \Pr{Y \leq m, ~X + Y > n} \\
			& + \Pr{X > m, Y > m, X + Y > n}.
\end{align*}
These three terms can be separately bounded from \emph{above} by $\Pr{Y > n-m}$, $\Pr{X > n-m}$ and $\Pr{X > m} \cdot \Pr{Y > n/2} + \Pr{Y > m} \cdot \Pr{X > n/2}$, respectively.
Again using $\lim_n \frac{a_n}{a_{n+m}} = 1$, we have that
\begin{align*}
    &\lim_{n \to \infty} a_n \cdot \Pr{Y > n-m} = \psi_a(\nu), \\
    &\lim_{n \to \infty} a_n \cdot \Pr{X > n-m} = \psi_a(\mu).
\end{align*}
And since $z := \sup_n \frac{a_n}{a_{\lfloor n/2\rfloor}} <\infty$, we get
\begin{align*}
    &\limsup_{n \to \infty} a_n \cdot \Pr{X > m} \cdot \Pr{Y > n/2} \leq z\epsilon \cdot \psi_a(\nu), \\
    &\limsup_{n \to \infty} a_n \cdot \Pr{Y > m} \cdot \Pr{X > n/2} \leq z\epsilon \cdot \psi_a(\mu). 
\end{align*}
Thus we obtain
\[
\limsup_{n \to \infty} a_n \cdot \eta((n,\infty)) \leq (1+z\epsilon)(\psi_a(\mu) + \psi_a(\nu)). 
\]
Letting $\epsilon \to 0$ then yields $\limsup_{n \to \infty} a_n \cdot \eta((n,\infty)) \leq \psi_a(\mu) + \psi_a(\nu)$. So $\psi_a(\eta) = \lim_{n \to \infty} a_n \cdot \eta((n,\infty))$ exists and equals the sum $\psi_a(\mu) + \psi_a(\nu)$. This shows $\psi$ is additive. 

\section{Other intermediate semigroups}

In this section, we exemplify the diverse behavior of the monotone additive statistics on various other sub-semigroups of $\P$.

\subsection{Extensions of the expectation}

We have shown that the expectation is the unique monotone additive statistic on $\P^1$, and that there are no monotone additive statistics on some semigroups that are larger than $\P^1$, including $\P$. In this subsection, we discuss some natural (and classical) extensions of the expectation to super-semigroups of $\P^1$.

Let $\P_\sigma$ be the set of probability measures $\mu$ on $\R$ such that
\begin{align*}
    \sigma(\mu) = \int_0^\infty [\mu((x,\infty))-\mu((-\infty,-x))]\,\dd x
\end{align*}
exists and is finite. Note that $\sigma(\mu) = \mathbb{E}[\mu]$ whenever $\mu \in \P^1$, but that $\P_\sigma$ is much larger than $\P^1$, and includes for example, all the symmetric measures. 

\begin{claim}
    $\P_\sigma$ is a semigroup, and $\sigma$ is a monotone additive statistic on it.
\end{claim}
\begin{proof}
    We first show that $\P_\sigma$ is a semigroup. Given $\mu \in \P_\sigma$,  denote by $\check \mu$ the reflection of $\mu$ at $0$, given by $\check\mu(A) := \mu(-A)$. Denote by $\mu_{\mathrm{even}}$ and $\mu_{\mathrm{odd}}$ the symmetric and anti-symmetric parts of $\mu$:
    \begin{align*}
        \mu_{\mathrm{even}} = \frac{\mu + \check \mu}{2}, \qquad \mu_{\mathrm{odd}} = \frac{\mu - \check \mu}{2}. 
    \end{align*}
	    Note that $\mu_{\mathrm{even}}$ is a symmetric probability measure on $\R$, and that $\mu_{\mathrm{odd}}$ is an anti-symmetric with total mass zero and total variation norm $\le 1$.

	    Define $A_\mu \colon \R \to \R$ by
	    \begin{align*}
		A_\mu(x) = \mu_{\mathrm{odd}}((x,\infty)) = \frac{\mu((x,\infty))-\mu((-\infty,-x))}{2}.
	    \end{align*}
	    Then 
	    \begin{align*}
		\sigma(\mu) = \int_{\R} A_\mu
	    \end{align*}
	    and $\mu \in \P_\sigma$ if and only if $A_\mu \in L^1(\R)$.

	    Suppose $\mu,\nu \in \P_\sigma$. We need to show that $\mu \ast \nu \in \P_\sigma$ and $\sigma(\mu \ast \nu) = \sigma(\mu) + \sigma(\nu)$. Note that
    \begin{align*}
        (\mu * \nu)_{\mathrm{odd}} = \frac{\mu * \nu - \check \mu * \check\nu}{2} = \mu_{\mathrm{even}} \ast \nu_{\mathrm{odd}} + \mu_{\mathrm{odd}} \ast \nu_{\mathrm{even}},
    \end{align*}
    and so
    \begin{align*}
        A_{\mu \ast \nu} = \mu_{\mathrm{even}} \ast A_\nu + \nu_{\mathrm{even}} \ast A_\mu.
    \end{align*}
    Since $\mu_{\mathrm{even}}$ and $\nu_{\mathrm{even}}$ are probability measures, $\mu_{\mathrm{even}} \ast A_\nu$ and $\nu_{\mathrm{even}} \ast A_\mu$ are in $L^1(\R)$, and furthermore $\int \mu_{\mathrm{even}} \ast A_\nu = \int A_\nu$ and likewise $\int \nu_{\mathrm{even}} \ast A_\mu = \int A_\mu$. Hence
    \begin{align*}
	    \sigma(\mu \ast \nu) = \int \left( \mu_{\mathrm{even}} \ast A_\nu + \nu_{\mathrm{even}} \ast A_\mu \right) = \int A_\mu + \int A_\nu = \sigma(\mu) + \sigma(\nu).
    \end{align*}
    Thus $\P_\sigma$ is a semigroup and $\sigma$ is additive. 
    The monotonicity of $\sigma$ is immediate from the definition.
\end{proof}

The monotone additive statistic $\sigma$ is related to Kolmogorov's generalized mathematical expectation \cite{kolmogorov2018foundations}*{p.\ 64}, sometimes called the weak mean, given by
\begin{align*}
    \mathbb{E}^*[\mu] = \lim_{n \to \infty}\int_0^n [\mu((x,\infty))-\mu((-\infty,-x))]\,\dd x
\end{align*}
whenever this limit exists and when $\lim_n n\cdot\mu\big((-\infty,n)\cup (n,\infty)\big) = 0$;  denote by $\P_*$ the set of probability measures $\mu$ satisfying these two conditions.

The significance of $\mathbb{E}^*$ is that it plays for the weak law of large numbers the role that the expectation plays for the strong law ~\cite{kolmogorov2018foundations}*{Section~VI.4}: Given i.i.d.\ random variables $X_1, X_2, \ldots$ with distribution $\mu$, the sequence $((X_1+\cdots+X_n)/n)_{n\ge 1}$ converges in probability to some constant $M$ if and only if $\mu \in \P_*$, in which case $M = \mathbb{E}^*[\mu]$. It follows immediately from this theorem that $\P_*$ is a semigroup, and that $\mathbb{E}^*$ is an additive statistic.\footnote{In this context, Pitman's theorem \cite{pitman1956derivatives} is also worth noting: $\mu \in \P_*$ if and only if the characteristic function $t \mapsto \int \ee^{i t x}\,\dd\mu(x)$ is differentiable at zero, in which case the derivative at zero is equal to $i\mathbb{E}^*[\mu]$.}
It is also straightforward to show that it is monotone.

Note that neither semigroup $\P_*$ and $\P_\sigma$ contains the other. The latter includes every symmetric distribution, including ones for which the second defining condition of $\P_*$ does not hold.
The former includes for example any measure $\mu$ with tails given by
\[
	\mu((-\infty,-x)) = \frac{1+\sin x}{x\log x}, \qquad \mu((x, \infty)) = \frac{1}{x\log x}
\]
for all $x \gg 1$.
Such measures are not in $\P_\sigma$, because $x \mapsto \mu((x,\infty))-\mu((-\infty,-x))$ is not integrable.

\subsection{Further examples}
\label{sec:further}
Let
\begin{align*}
    \psi(\mu) = \lim_n n \mu((n,\infty)).
\end{align*}
By Theorem~\ref{thm:Lpsi}, this is a monotone additive statistic on the semigroup of probability measures on which it is finite.

Given a measure $\mu$, denote by $\check \mu$ the reflection of $\mu$ at $0$, given by $\check\mu(A) := \mu(-A)$. Let $\P^{\psi}_{\pm}$ be the semigroup of measures $\mu$ for which both $\psi(\mu)$ and $\psi(\check \mu)$ exist and are finite. This semigroup lies strictly between $\P^1$ and $\P^p$ for any $p < 1$. Despite the fact that zero is the only monotone additive statistic on these $\P^p$ semigroups, for any $a,b \ge 0$ the assignment
\[
	\mu \longmapsto a\psi(\mu) - b\psi(\check\mu)
\]
defines a monotone additive statistic on $\P^{\psi}_{\pm}$. We do not know if these are the only monotone additive statistics on this semigroup.

We can obtain other interesting examples by considering even smaller semigroups that still contain $\P^1$. Let $\P^\psi_=$ be the set of $\mu \in \P^\psi_\pm$ for which $\psi(\mu) = \psi(\check \mu)$. This is a semigroup strictly between $\P^{1}$ and $\P^{\psi}_{\pm}$. Moreover, since $\mu \geq \nu$  only if $\psi(\mu) \geq \psi(\nu)$ and $\psi(\check \mu) \leq \psi(\check \nu)$, stochastic dominance on $\P^\psi_=$ requires $\psi(\mu) = \psi(\check\mu) = \psi(\nu) = \psi(\check\nu)$. Now take $\gamma : \R \to \R$ to be any additive function, i.e., any function that satisfies the Cauchy functional equation $\gamma(x+y)=\gamma(x)+\gamma(y)$. Then $\gamma \circ \psi$ is an additive statistic on $\P^\psi_=$, and is also trivially monotone. 
Thus the set of monotone additive statistics on the semigroup $\P^\psi_=$ is rather complex.

We can consider the smaller semigroup $\P^\psi_0 \subset \P^\psi$ of measures $\mu$ for which $\psi(\mu) = \psi(\check\mu) = 0$. We claim that zero is the only monotone additive statistic on this semigroup; the proof is almost identical to the proof of Theorem~\ref{thm:L1}, except that the $a_n$ are chosen large enough so that $x \cdot \nu_n((x,\infty)) \leq \frac{1}{n}$ for all $x \geq a_n$.

\subsection{Overview of semigroups considered}

The inclusion relationships among most of the semigroups discussed in this paper can be nicely summarized as follows. Here $q$ is assumed to be in $(1,\infty)$, and $p \in (0,1)$.
\begin{align*}
	\P_c \subset \P_{\mathrm{Cram}} \subset \bigcap_{q > 1} \P^q \subset & \; \P^1 \\ & \subset  \P^\psi_0 \subset \P^\psi_= \subset \P^\psi_{\pm} \subset \P^{p} \subset \P,
\end{align*}
Note that all inclusions are strict.
Additionally, the semigroup $\P^\psi$ is strictly intermediate between $\P^1$ and $\P$.

For those semigroups listed here that strictly contain $\P^1$ (second row), we have shown that $\P^\psi_0$ and $\P^{p}$ for any $p \in (0,1)$, as well as the limit case $\P^0 = \P$, only admit the zero monotone additive statistic.
However, the intermediate semigroups $\P^\psi_=$ and $\P^\psi_{\pm}$ have non-trivial monotone additive statistics, and the former even has a very large set thereof.

\section{Open questions}
We end the paper with a number of open questions. 
\begin{itemize}
    \item The variance is an additive statistic which instead of monotonicity has the related property of non-negativity. Is there a (non-trivial) non-negative additive statistic on the entirety of $\P$? 
    \item While the variance is not monotone with respect to the stochastic order, it is monotone with respect to the convex order. Is it the only additive statistic on $\P^2$ (and on smaller semigroups) which is monotone for the convex order? \item In this paper we study convolutions. A related natural operation is the averaging of independent random variables, which has produced interesting recent results involving heavy tailed distributions \cite{chen2025unexpected}. Are there similar questions that can be asked about this operation?
    \item Mattner~\cite{mattner2004cumulants} shows that the variance is the unique continuous non-negative additive statistic on $\bigcap_p \P^p$. Are there non-continuous ones?
    \item The stochastic order can be defined for $\R^d$ with any $d \in \N$, equipped with the natural (product) partial order. What are the monotone additive statistics on $\P(\R^d)$?
    \item Given $\mu,\nu \in \P$, under what conditions is there $\eta \in \P$ with $\mu \ast \eta \ge \nu \ast \eta$? For $\mu, \nu \in \P^1$, it was shown in \cite{PST} that if $\mu$ and $\nu$ are different, then a necessary and sufficient condition is that $\E{\mu}>\E{\nu}$.\footnote{Sufficiency is Theorem 1 in \cite{PST}; the end of their Appendix A shows necessity.}
    \item What are the additive statistics of $\P^{\psi}_{\pm}$, as defined in Section \ref{sec:further}?
    \item For the sub-semigroup $\P(\R_+)$ of probability measures on $[0,\infty)$, the expectation is a monotone additive statistic taking values in $[0,\infty]$. This is not true for $\P$, where the expectation is ill-defined even if $\pm \infty$ is allowed. A natural question is to understand the monotone additive statistics $\P(\R_+) \to [0,\infty]$. This is a rich collection which include many of the examples presented in this paper.
    
\end{itemize}

\bibliography{refs}

\end{document}